\documentclass[10pt]{amsart}
\usepackage{amsmath,amssymb,amsrefs}

\numberwithin{equation}{section}

\newtheorem{theorem}[equation]{Theorem}
\newtheorem{lemma}[equation]{Lemma}
\newtheorem{corollary}[equation]{Corollary}
\newtheorem{prop}[equation]{Proposition}
\newtheorem{question}[equation]{Question}

\theoremstyle{definition}
\newtheorem{definition}[equation]{Definition}

\theoremstyle{remark}
\newtheorem{remark}[equation]{Remark}

\begin{document}

\title{Capability of nilpotent products of cyclic groups}
\author{Arturo Magidin}

\subjclass[2000]{Primary 20D15, Secondary 20F12}

\maketitle

\begin{center}
\textit{Department of Mathematics}\\
\textit{University of Louisiana at Lafayette}\\
\textit{217 Maxim Doucet Hall}\\
\textit{P.O. Box 41010}\\ 
\texttt{magidin@member.ams.org}
\end{center}

\bigskip
Published in \textsl{J.~Group Theory} {\bf 8} (2005), pp.~431--452.

\bigskip
Errata added on 4/10/2006 -- Updated on 5/6/2006

\vfill\eject

{\small
\noindent\textit{Abstract.} A group is called capable if it is a central factor
group. We consider the capability of nilpotent products of cyclic
groups, and obtain a generalization of a theorem of Baer for the small
class case. The approach
is also used to obtain some recent results on the capability
of certain nilpotent groups of class~$2$. We also prove a necessary
condition for the capability of an arbitrary $p$-group of class $k$,
and some further results.
}

\section{Introduction}

In his landmark paper on the classification of finite
$p$-groups~\cite{hallpgroups}, P.~Hall remarked that:
\begin{quote}
The question of what conditions a group $G$ must fulfil in order that
it may be the central quotient group of another group $H$, $G\cong
H/Z(H)$, is an interesting one. But while it is easy to write down a
number of necessary conditions, it is not so easy to be sure that they
are sufficient.
\end{quote}
Following M.~Hall and Senior~\cite{hallsenior}, we make the following
definition:

\begin{definition} A group $G$ is said to be \textit{capable} if and only if
there exists a group $H$ such that $G\cong H/Z(H)$; equivalently, if
and only if $G$ is isomorphic to the inner automorphism group of a
group~$H$.
\end{definition}

Capability of groups was first studied by R.~Baer in~\cite{baer},
where, as a corollary of some deeper investigations, he characterised
the capable groups that are a direct sum of cyclic groups.  Capability
of groups has received renewed attention in recent years, thanks to results of
Beyl, Felgner, and Schmid~\cite{beyl} characterising the capability of
a group in terms of its epicenter; and more recently to work of Graham
Ellis~\cite{ellis} that describes the epicenter in terms of the
nonabelian tensor square of the group. The epicenter was used
in~\cite{beyl} to characterize the capable extra-special $p$-groups;
and the nonabelian tensor square was used in~\cite{baconkappe} to
characterize the capable $2$-generator finite $p$-groups of odd order
and class~$2$.

While the nonabelian tensor product has proven very useful in the
study of capable groups, it does seem to present certain
limitations. At the end of~\cite{baconkappe}, for example, the authors
note that their methods require ``very explicit knowledge of the
groups'' in question.

In the present work, we will use ``low-tech'' methods to obtain a
number of results on capability of finite $p$-groups. They rely only
on commutator calculus, and so may be more susceptible to
extension than results that require explicit knowledge of the
nonabelian tensor square of a group. In particular, we will prove a
generalization of Baer's Theorem characterising the capable direct
sums of cyclic groups to the $k$-nilpotent products of cyclic
$p$-groups with $k<p$, and with $k=p=2$. 

One weakness in our results should be noted explicitly: Baer's Theorem
fully characterizes the capable finitely generated abelian groups,
because every finitely generated abelian group can be expressed as a
direct sum of cyclic groups. Our generalization does not provide a
result of similar reach, because whenever $k>1$ there exist finite
$p$-groups of class~$k$ that are not a $k$-nilpotent product of cyclic
$p$-groups, even if we restrict to $p>k$. On the other hand, every
finite $p$-group of class $k$ is a \textit{quotient} of a
$k$-nilpotent product of cyclic $p$-groups, so some progress can be
made from this starting point.

The paper is organized as follows: in Section~\ref{sec:defs} we
present the basic definitions and notation. We proceed in
Section~\ref{sec:nec} to establish a necessary condition for
capability which extends an observation of P.~Hall. In
Section~\ref{sec:centkprod}, we first describe the center of a
$k$-nilpotent product of cyclic $p$-groups, where $p$ is a prime
satisfying $p\geq k$, and then use this description to prove
Theorem~\ref{capabilitynilkprod}, the promised generalization of
Baer's Theorem.  In
Section~\ref{sec:kp2} we characterize the capable $2$-nilpotent
products of cyclic $2$-groups. Then in 
Section~\ref{sec:applic} we use our results on capable 2-nilpotent
products of cyclic $p$-groups to derive a characterization of the
capable $2$-generated nilpotent $p$-groups of class two, $p$ an odd
prime (recently obtained through different methods by Bacon
and Kappe), and some other related results, 
by way of illustration of how our methods may be
used as a starting point for further investigations. 

The main results are Theorem~\ref{necessity}, giving a necessary
condition for capability of a finite $p$ group of class~$k$, and
Theorems \ref{capabilitynilkprod} and~\ref{pandkequal2}, characterising the capable
$k$-nilpotent products of cyclic $p$-groups for $p>k$ and for~$k=p=2$.

\section{Definitions and notation\label{sec:defs}}

All maps will be assumed to be group homomorphisms. All groups will be
written multiplicatively, unless we explicitly state otherwise.

Let $G$ be a group. The center of~$G$ is denoted by~$Z(G)$.
The identity element of~$G$ will be denoted by $e$.

Let $x\in G$. We say that $x$ is \textit{of exponent
  $n$} if and only if $x^n=e$; we say $x$ is \textit{of order $n>0$}
  if and only if $x^n=e$ and $x^k\neq e$ for all $k$, $0<k<n$. 

The commutator $[x,y]$ of two elements $x$ and $y$ is
defined to be $[x,y]=x^{-1}y^{-1}xy$; given two subsets (not necessarily
subgroups) $A$ and~$B$ of~$G$, we let $[A,B]$ be the subgroup
generated by all elements of the form $[a,b]$, with $a\in A$ and $b\in
B$.

The lower central series of~$G$ is the sequence of subgroups defined
by $G_1,G_2,\ldots$, where $G_1=G$, $G_{n+1}=[G_n,G]$. We say a group
$G$ is \textit{nilpotent of class (at most)~$n$} if and only if
$G_{n+1}=\{e\}$; we will often drop the \textit{``at most''} clause, it
being understood. The nilpotent groups of class~$1$ are the abelian
groups.  Note that $G$ is nilpotent of class~$n$ if and only if $G_n\subset
Z(G)$. The class of all nilpotent groups of class at most~$k$ is denoted
by~$\germ N_k$; it is a variety of groups in the sense of General
Algebra; we refer the reader to Hanna Neumann's excellent
book~\cite{hneumann}.

We will write commutators left-normed, so that
$[a_1,a_2,a_3] = [[a_1,a_2],a_3]$,
etc. 

The following properties of the lower central series are
well-known. See for example~\cite{hall}:

\begin{definition} Let $g\in G$, $g\not= e$. We define $W(g)$ to be
  $W(g)=k$ if and only if $g\in G_k$ and $g\not\in G_{k+1}$. We also
  set $W(e)=\infty$.
\end{definition}

\begin{prop} Let $G$ be a group.
\begin{itemize}
\item[(i)] For all $a,b\in G$, $W([a,b])\geq W(a) + W(b)$.
\item[(ii)] If $W(a_i)=w_1$ and $W(b_j)=w_2$, then
\[ \left[ \prod_{i=1}^I a_i^{\alpha_i}\, , \,\prod_{j=1}^J
  b_j^{\beta_j}\right] \equiv \prod_{i=1}^I\prod_{j=1}^J
  [a_i,b_j]^{\alpha_i\beta_j} \pmod {G_{w_1+w_2+1}}.\]
\item[(iii)] If $a\equiv c \pmod{G_{W(a)+1}}$ and $b\equiv d
  \pmod{G_{W(b)+1}}$, then
\[ [a,b] \equiv [c,d] \pmod{G_{W(a)+W(b)+1}}.\]
\item[(iv)] A variant of the Jacobi identity:
\[ [a,b,c]\,[b,c,a]\,[c,a,b] \equiv e \pmod{G_{W(a)+W(b)+W(c)+1}}.\]
\end{itemize}
\label{Widentities}
\end{prop}

The following identities may be verified by direct calculation:
\begin{eqnarray}
\null[xy,z] & = & [x,z][x,z,y][y,z]\label{prodformone}\\
\null[x,yz] & = & [x,z][z,[y,x]][x,y].\label{prodformtwo}\\
\null[a^r,b^s] & \equiv & [a,b]^{rs} [a,b,a]^{s
\binom{r}{2}}
[a,b,b]^{r\binom{s}{2}}
\pmod{G_4}\label{prodformthree}\\
\null[b^r,a^s] &\equiv& [a,b]^{-rs} [a,b,a]^{-r\binom{s}{2}}[a,b,b]^{-s\binom{r}{2}}
\pmod{G_4},\label{prodformfour}
\end{eqnarray}
where $\binom{r}{2} = \frac{r(r-1)}{2}$ for all integers $r$.

\subsection*{Nilpotent product of groups}

The nilpotent products of groups were introduced by
Golovin~\cite{golovinnilprods} as examples of regular products of
groups. Although defined in a more general context in which there is
no restriction on the groups involved, our definition will be
restricted to the situation we are interested in.

\begin{definition} Let $A_1,\ldots,A_n\in{\germ N}_k$. The
\textit{$k$-nilpotent product of $A_1,\ldots,A_n$,} denoted by
$A_1 \amalg^{\germ N_k} \cdots \amalg^{\germ N_k} A_n$,
is defined to be the group $G=F/F_{k+1}$, where $F$ is the free
product of the $A_i$, $F=A_1 * \cdots * A_n$, and $F_{k+1}$ is the
$(k+1)$-st term of the lower central series of $F$.
\end{definition}

Note that the ``$1$-nilpotent product'' is simply the direct
sum. Also, if the $A_i$ are in $\germ N_{k-1}$, and $G$ is the
$k$-nilpotent product of the $A_i$, then the $(k-1)$-nilpotent product
of the $A_i$ is isomorphic to $G/G_k$.

The use of the coproduct notation does not appear to be standard in
the literature, but there is a good reason to use it: the $k$-nilpotent
product as defined above is the coproduct (in the sense of category
theory) in the category $\germ N_{k}$.

\subsection*{Basic commutators}

The collection process of M.~Hall gives normal forms for relatively
free groups in ${\germ N_k}$, and in some other special cases. The
concept of basic commutators is essential in this developement. The
definition of basic commutators seems to vary in the
literature (and sometimes even within the same work; cf. \S 11.1 and
\S 12.3 in~\cite{hall}). The definition we will use in the
present work is that which appears in \S 12.3 of~\cite{hall}, where
the ordering of commutators of weight~$n>1$ is given by letting
$[x_1,y_1]<[x_2,y_2]$ if and only if $y_1<y_2$ or $y_1=y_2$ and
$x_1<x_2$ (lexicographically from right to left); here $[x_1,y_1]$
and~$[x_2,y_2]$ are basic commutators of weight~$n$.

For an exposition of the collection process, we direct the reader to
Chapter~11 in~\cite{hall}. The main consequences are Hall's Basis
Theorem~\cite{hall}*{Theorem 11.2.4} and the collection
formulas~\cite{hall}*{\S 12.3}. We will use a slightly more precise
version of the collection formulas in Section~\ref{sec:nec} below. We
will also use later on a special case of the generalization of the
Basis Theorem due to R.R.~Struik \cite{struikone}*{Theorem~3}; it
states that if we take the $k$-nilpotent product of $n$
cyclic~$p$-groups, $p\geq k$, generated by $x_1,\ldots,x_n$, then
there is a normal form for the elements as products $\prod c_i^{a_i}$,
where $c_1<c_2<\cdots$ are the basic commutators on the $x_i$ of
weight less than or equal to~$k$, and the exponent~$a_i$ is taken
modulo the smallest of the orders of the generators that appear in the
full expression of~$c_i$.

A corollary of this result is that if $F$ is the $k$-nilpotent product
of the cyclic~$p$-groups, $p\geq k$, then $F_i/F_{i+1}$ is abelian
with basis given by the basic commutators of weight exactly $i$,
$i=1,\ldots,k$. 

\section{A necessary condition\label{sec:nec}}

In this section we give a necessary condition for capability of finite
$p$-groups based on the orders of the elements on a minimal generating
set.  In the case of small class (that is, when $G$ is a $p$-group
with $G\in {\germ N}_k$ and $p>k$), the condition reduces to an observation which
goes back at least to P.~Hall (penultimate paragraph in pp.~137
in~\cite{hallpgroups}). Although Hall only considers bases in the
sense of his theory of regular $p$-groups, his argument is essentially
the same as the one we present.  However, Hall's result may not be
very well known, since it is only mentioned in passing; see for
example~\cite{baconkappe}*{Theorem~4.4}.

We begin by recalling three consequences of the collection process:

\begin{lemma}[Lemma H1 in~\cite{struikone}]
Let $x,y$ be any elements of a group; let $c_1,c_2,\ldots$ be the
sequence of basic commutators of weight at least two in $x$ and~$[x,y]$,
in ascending order. Then
\begin{equation}
[x^{\alpha},y] = [x,y]^{\alpha}c_1^{f_1(\alpha)} c_2^{f_2(\alpha)}\cdots
c_i^{f_i(\alpha)}\cdots
\label{approxformula}
\end{equation}
where
\begin{equation}
 f_i(\alpha) = a_1\binom{\alpha}{1} + a_2\binom{\alpha}{2} +
\cdots + a_{w_i}\binom{\alpha}{w_i},
\label{formofthefis}
\end{equation}
where $a_i\in\mathbb{Z}$,
and $w_i$ is the weight of $c_i$ in $x$ and $[x,y]$.  If the group is
nilpotent, then the expression in {\rm(\ref{approxformula})} gives an
identity, and the sequence of commutators terminates.
\label{struiklemmah1}
\end{lemma}

\begin{lemma} [Lemma H2 in~\cite{struikone}]
Let $\alpha$ be a fixed integer, and $G$ a nilpotent group of class at
most $k$. If $b_j\in G$ and $r<k$, then
\begin{equation}
[b_1,\ldots,b_{i-1},b_i^{\alpha},b_{i+1},\ldots,b_r] =
   [b_1,\ldots,b_r]^{\alpha} v_1^{f_1(\alpha)}
   v_2^{f_2(\alpha)}\cdots v_t^{f_t(\alpha)}
\label{pullingoutexp}
\end{equation}
 where the $v_k$ are (not necessarily basic) commutators in $b_1,\ldots,b_r$ of weight
   strictly greater than $r$, and every $b_j$, $1\leq j\leq r$ appears
   in each commutator $v_k$.
   Furthermore, the $f_i$ are of the form~{\rm (\ref{formofthefis})},
   with $a_j\in\mathbb{Z}$, and $w_i = w_i' - (r-1)$, where $w_i'$ is
   the weight of $v_i$ in the $b_i$.
\label{struiklemmah2}
\end{lemma}

We find {\rm(\ref{pullingoutexp})} useful in situations when we have
commutators in some terms, some of which are shown as powers, and we
want to ``pull the exponent out.'' At other times, we want to
reverse the process and pull the exponents ``into'' a commutator. In
such situations, we use {\rm ({\ref{pullingoutexp})}} to express
$[b_1,\ldots,b_r]^{\alpha}$ in terms of other commutators. We will call the
resulting identity $(\mathrm{\ref{pullingoutexp}}')$; that is:
\begin{equation}
\relax[b_1,\ldots,b_r]^{\alpha} =
      [b_1,\ldots,b_{i-1},b_i^{\alpha},b_{i+1},\ldots,b_r]v_t^{-f_t(\alpha)}\cdots
      v_2^{-f_2(\alpha)}v_1^{-f_1(\alpha)},\tag{$\mathrm{\ref{pullingoutexp}}'$}
\end{equation}
with the understanding that we will only do this in a nilpotent group
so that the formula makes sense. 

\begin{lemma}[Theorem 12.3.1 in~\cite{hall}]
Let $x_1,\ldots,x_s$ be any $s$ elements of a group. Let
$c_1,c_2,\ldots$ be the basic commutators in $x_1,\ldots,x_s$ of
weight at least~$2$, written in increasing order. Then
\begin{equation}
(x_1\cdots x_s)^{\alpha} = x_{1}^{\alpha} x_{2}^{\alpha}\cdots x_{s}^{\alpha}
c_1^{f_1(\alpha)}\cdots c_i^{f_i(\alpha)}\cdots
\label{approxforpower}
\end{equation}
where $f_i(\alpha)$ is of the form~{\rm(\ref{formofthefis})}, with 
$a_j$ integers that depend only on~$c_i$ and not on~$\alpha$, and $w_i$ is the
weight of $c_i$ in the $x_j$.  If the group is nilpotent, then
equation {\rm(\ref{approxforpower})} is satisfied as an identity in the group,
and the sequence of commutators terminates. 
\label{struiktheoremh3}
\end{lemma}

The following lemma is easily established by induction on the weight:

\begin{lemma} Let $F(y,z)$ be the free
  group on two generators. Then every basic commutator of ${\rm
  weight}\geq 3$ is of the form
\begin{equation}
\relax[z,y,y,c_4,\ldots,c_r]\quad\mbox{or}\quad[z,y,z,c_4,\ldots,c_r]
\label{formofbasic}
\end{equation}
where $r\geq 3$, and $c_4,\ldots,c_r$ are basic commutators in $y$
and~$z$.
\label{descriptionbasic}
\end{lemma}

The main idea in our development is as follows: if we know that
$[z,y^{p^i}]$ centralizes $\langle y,z\rangle$ in a group $G\in{\germ
N}_k$, for some prime $p$ and all integers~$i$ greater than or equal
to a given bound~$a$, then we want to prove that a commutator of the
form $[z^{p^n},y]$ is equal to $[z,y^{p^n}]$. To accomplish this, we
observe that a basic commutator of weight $k$ will have exponent
$p^a$, since we may simply use Lemma~\ref{descriptionbasic} and
Proposition~\ref{Widentities}(ii) to pull the exponent into the second
entry of the bracket. An arbitrary commutator of weight $k$ will also
have the same exponent, since $G_k$ is abelian. For a basic commutator
of weight $k-1$ we may use $(\mathrm{\ref{pullingoutexp}}')$ and
deduce that a sufficiently high power of $p$ will again yield the
trivial element, by bounding below the power of $p$ that divides the
exponents $f_i(p^n)$.  Then we apply Lemma~\ref{struiktheoremh3} to
deal with an arbitrary element of $G_{k-1}$.  Continuing in this way,
we can show that $\langle y,z\rangle_3$ is of exponent $p^M$ for
some~$M$, and we will obtain the desired result by applying
Lemma~\ref{struiklemmah1} to $[z^{p^N},y]$ and $[z,y^{p^N}]$ for a
suitably chosen large~$N$, thus showing that they are both equal to
$[z,y]^{p^N}$. Most of the work will go into obtaining a good estimate
on how large this ``large $N$'' has to be for everything to work.

If $p$ is a prime and $a$ is a positive integer, 
we let $[a]_p$ denote the exact $p$-divisor of~$a$; that is, we say
that $[a]_p=r$ if $p^r|a$ and $p^{r+1}\not|a$. Formally, we set $[0]_p
= \infty$. A classical theorem of Kummer implies that if $a$ is a
positive integer, $0<a\leq p^n$, then
\[ \left[ \binom{p^n}{a}\right]_p = n - [a]_p.\]
Recall that if $x$ is a real number, then $\lfloor x\rfloor$
  denotes the floor of~$x$, the largest integer smaller than or equal to~$x$.

\begin{prop} Let $p$ be a prime, $n$ a positive integer, and $m$
  an integer with $0<m\leq p^n$. If $a_1,\ldots,a_m$ are integers,
  then
\begin{equation}
 a_1\binom{p^n}{1} + a_2\binom{p^n}{2} + \cdots +
  a_m\binom{p^n}{m}\equiv 0 \pmod{p^{n-d}},
\label{sumofbinoms}
\end{equation}
where $d$ is the smallest integer such that
$p^{d+1}>m$; that is, $d=\lfloor\log_p(m)\rfloor$.
\label{boundfis}
\end{prop}

\begin{proof} Write $m=m_0 + m_1p + \cdots + m_dp^d$, with $0\leq
  m_i<p$ and $m_d>0$. Then for all integers $k$ between $1$ and $m$,
  $0\leq [k]_p\leq d$. Therefore, $n-[k]_p\geq n-d$,
  so each summand in {\rm(\ref{sumofbinoms})} is divisible by $p^{n-d}$,
  as claimed. 
\end{proof}

\begin{lemma}
Let $k\geq 3$, and let $G\in{\germ N}_k$. Let $p$ be a prime, $a\geq 0$ an
integer, and $y,z\in G$. Suppose that
\begin{equation}
\forall i\geq a,\qquad [z,y^{p^i},y] = [z,y^{p^i},z] = e.\label{identityppower}
\end{equation}
Then $\langle y,z\rangle_k$ is of exponent $p^a$, and $\langle
y,z\rangle_{k-1}$ is of exponent $p^{a+\left\lfloor\frac{1}{p-1}\right\rfloor}$.
\label{inductionbase}
\end{lemma}

\begin{proof}
It follows by $(\mathrm{\ref{pullingoutexp}}')$, 
(\ref{formofbasic}), and our assumption that
\[\relax [z,y,w,c_4,\ldots]^{p^a} = [z,y^{p^a},w,c_4,\ldots] = e\]
where $w\in \{y,z\}$, and $[z,y,w,c_4,\ldots]$ is a basic commutator of weight~$k$.
Since $\langle y,z\rangle_k$ is abelian and generated by the basic
commutators, this proves the first part of the statement.

Now let $c\in\langle y,z\rangle_{k-1}$ be a basic commutator of weight
exactly $k-1$, and let $N\geq a$. Expressing $c$ in the form given in
(\ref{formofbasic}), say $c = [z,y,w,c_r,\ldots]$ with $w\in\{y,z\}$, then
$(\mathrm{\ref{pullingoutexp}}')$ and our assumption yield
\[ c^{p^N} = [z,y,w,c_4,\ldots]^{p^N} =
   v_t^{-f_t(p^N)}\cdots v_1^{-f_1(p^N)},\]
where all the $v_i$ are of weight $k$. 
The corresponding exponents are of the form 
\[f_i(p^N)= a_1\binom{p^N}{1} + a_2\binom{p^N}{2}.\]
By Proposition~\ref{boundfis}, this expression is divisible by
$p^{N-\lfloor \log_p(2)\rfloor}$, and so whenever $N\geq a+\lfloor
\log_p(2)\rfloor$, we have $c^{p^N}=e$. Therefore, every basic
commutator of weight exactly $k-1$ is of exponent $p^{a+\lfloor
\log_p(2)\rfloor}$. Since $k>2$, $G_{k-1}$ is abelian as well; as every
generator is of exponent $p^{a+\lfloor \log_p(2)\rfloor}$, the theorem
is now proven by noting that $\lfloor\log_p(2)\rfloor =
\left\lfloor\frac{1}{p-1}\right\rfloor$.
\end{proof}

We will need the following technical lemma:

\begin{lemma} Let $k$ and $n$ be positive integers, with $n>1$. Then
\[ \max_{1\leq s\leq k} \left( \left\lfloor
\frac{k-s}{n-1}\right\rfloor +
\left\lfloor\log_n(s+1)\right\rfloor\right) = \left\lfloor
\frac{k}{n-1}\right\rfloor.\]
If $k\geq n-1$
then the maximum is always
attained at $s=n-1$.
\label{boundfors}
\end{lemma}

\begin{proof} Fix $k$ and $n$; we may assume $k\geq n-1$. Let 
\[ f(s) = \frac{k-s}{n-1} + \lfloor\log_n(s+1)\rfloor,\quad s\in [1,k].\]
The function $f(s)$ is piecewise strictly decreasing on $[1,k]$, and
continuous from the right at every point. So the maximum of $f(s)$
on $[1,k]$ will be achieved at $s=1$, or at a point where $f(s)$ is not
continuous. The points of discontinuity are the points of the form
$s=n^r-1$ where $1\leq r\leq \lfloor\log_n(k+1)\rfloor$. At $s=1$, we
have 
\[ f(1)=\frac{k-1}{n-1}+\lfloor\log_n(2)\rfloor \leq \frac{k}{n-1},\]
with equality if and only if $n=2$. At the points of discontinuity we
have
\begin{eqnarray*}
f(n^r-1) &=& \frac{k}{n-1} - \frac{n^r-1}{n-1} + r \\
&=& \frac{k}{n-1} - (1+n+n^2+\cdots+n^{r-1}) + r
\leq \frac{k}{n-1},
\end{eqnarray*}
with equality if and only if $r=1$. Thus the maximum value of $f(s)$
is equal to $\frac{k}{n-1}$, always achieved at $s_0=n-1$. The result
now follows by taking $\lfloor f(s)\rfloor$.
\end{proof}

\begin{lemma}[cf.~Lemma 8.83 in~\cite{nildoms}] 
Let $G$ be a group, with $G\in{\germ N}_k$ and $k\geq 3$.
Furthermore, let $p$ be a prime,
$a>0$ an integer, and $y,z\in G$. Assume that $a$, $y$, $z$, and $G$
satisfy {\rm (\ref{identityppower})}.
Then $\langle y,z\rangle_{k-m}$ is of exponent $p^{a+\left\lfloor\frac{m}{p-1}\right\rfloor}$
for $m=0,1,\ldots,(k-3)$.
\label{expofthirdterm}
\end{lemma}

\begin{proof} We proceed by induction on
  $m$. Lemma~\ref{inductionbase} proves cases $m=0,1$, so we let
  $m\geq 2$.

Assume the result is true for $j=0,1,\ldots,m-1$. First, we consider a
basic commutator $c$ of weight exactly $k-m$. We may express $c$ as in
{\rm (\ref{formofbasic})}, and applying
$(\mathrm{\ref{pullingoutexp}}')$ to $c^{p^N}$ we obtain
\[ c^{p^N} = [z,y^{p^N},w,c_4,\ldots,c_r]v_t^{-f_t(p^N)}\cdots v_2^{-f_2(p^N)}
v_1^{-f_1(p^N)},\]
where $w\in\{y,z\}$, $c_4,\ldots,c_r$ are basic commutators in $z$ and
$y$, and the $v_i$ are commutators in $z,y,w,c_4,\ldots,c_r$;
we know that each of $z,y,w,c_4,\ldots,c_r$ appears at least once in
each $v_i$, and that we may express the weight of $v_i$ on
$z,y,w,c_4,\ldots,c_r$
as $r+s$ for some positive integer~$s$. Thus,
we may conclude that if the weight of $v_i$ in $z,y,w,c_4,\ldots,c_r$
is $r+s$, then $W(v_i)\geq (k-m)+s$. In particular, we consider only those
commutators with $1\leq s\leq m$. 

If $v_i$ is of weight $r+s$ in $z,y,w,c_4,\ldots,c_r$, then by
Lemma~\ref{struiklemmah2} we have that
\[ f_i(p^N) = a_1\binom{p^N}{1}+\cdots +a_{s+1}\binom{p^N}{s+1}.\]
Therefore, $f_i(p^N)$ is divisible by $p^{N-\lfloor
\log_p(s+1)\rfloor}$. Since $v_i\in\langle y,z\rangle_{(k-(m-s))}$, by
the induction hypothesis we know that $v_i^{-f_i(p^N)}$ is trivial whenever
$N-\lfloor\log_p(s+1)\rfloor \geq a+\left\lfloor\frac{m-s}{p-1}\right\rfloor$. Therefore, if
\[ N \geq a + \left\lfloor\frac{m-s}{p-1}\right\rfloor + \lfloor\log_p(s+1)\rfloor\]
for $s=1\ldots,m$, then we may conclude that $c^{p^N}=e$. By
Lemma~\ref{boundfors}, the largest of these values is
$a+\lfloor\frac{m}{p-1}\rfloor$, which shows that
the basic commutators of weight exactly $m-k$ are of
exponent $p^{a+\lfloor\frac{m}{p-1}\rfloor}$, as desired.

Now take an arbitrary element~$c$ of $\langle y,z\rangle_{(k-m)}$, and
write $c = d_1\cdots d_r$, where $d_i=c_i^{\beta_i}$ is the power of a
basic commutator of weight at least $k-m$ in $y$ and~$z$, and
$c_1<\cdots< c_r$. To calculate $c^{p^N}$, we apply
Lemma~\ref{struiktheoremh3} to this expression. We obtain:
\[(d_1\cdots d_r)^{p^N} = c_{1}^{\beta_1 p^N} c_{2}^{\beta_2
  p^N}\cdots c_{r}^{\beta_r p^N}
u_1^{f_1(p^N)}\cdots u_i^{f_i(p^N)}\cdots
\]
where $f_i(p^N)$ is of the form~{\rm(\ref{formofthefis})}, with $w_i$
the weight of $u_i$ in the $d_j$. If we let the weight of $u_i$ in the
$d_j$ be equal to $s$, then we know that $2\leq s$.  Since $W(d_i)\geq
(m-k)$, we have that $W(u_i)\geq s(m-k)$. In particular, we may
restrict our attention to values of $s$ satisfying $2\leq s\leq
\left\lfloor\frac{k}{k-m}\right\rfloor$.  

If $u_i$ is of weight $s$ in the $d_j$, then it lies in $\langle
y,z\rangle_{k-(m-(s-1)(k-m))}$, so by the induction hypothesis it is
of exponent $p^{a+\lfloor\frac{m-(s-1)(k-m)}{p-1}\rfloor}$. By (\ref{formofthefis})
and Proposition~\ref{boundfis} we have that $f_i(p^N)$ is a multiple
of $p^{N-\lfloor\log_p(s)\rfloor}$; thus we can guarantee that
$u_i^{f_i(p^N)}$ is trivial if
\[ N \geq a + \left\lfloor\frac{m-(s-1)(k-m)}{p-1}\right\rfloor + \lfloor\log_p(s)\rfloor.\]
Since $s\geq 2$ and $k-m\geq 3$, it is clear that this value will
certainly be no larger than $a+\lfloor\frac{m}{p-1}\rfloor$ by
Lemma~\ref{boundfors}, so we 
conclude that any element of $\langle y,z\rangle_{(k-m)}$ is of
exponent $p^{a+\lfloor\frac{m}{p-1}\rfloor}$, as claimed.
\end{proof}

\begin{theorem}[cf. Corollary~8.84 in~\cite{nildoms}]
Let $k\geq 1$ and $G\in{\germ N}_{k+1}$. Furthermore, let $p$ be a prime, $a>0$ an
integer, and $y$ and~$z$ elements of~$G$. Assume that
\[ \forall i\geq a,\qquad [z,y^{p^i},y]=[z,y^{p^i},z]=e.\]
If $N\geq a+\left\lfloor\frac{k-1}{p-1}\right\rfloor$, then
$[z^{p^N},y] = [z,y]^{p^N} = [z,y^{p^N}]$.
\label{likenil2withpowers}
\end{theorem}

\begin{proof} Applying Lemma~\ref{struiklemmah1} to $[z^{p^N},y]$ yields:
\[ [z^{p^N},y] = [z,y]^{p^N} v_1^{f_1(p^N)}v_2^{f_2(p^N)}\cdots v_t^{f_t(p^N)},\]
where $v_1,v_2,\ldots$ are the basic commutators of weight at least
$2$ in $z$ and $[z,y]$. If $v_i$ is of weight $s\geq 2$ in $z$ and
$[z,y]$, then 
\[ f_i(p^N) = a_1\binom{p^N}{1} + a_2\binom{p^N}{2} + \cdots +
a_s\binom{p^N}{s},\]
and we know that $v_i\in\langle y,z\rangle_{s+1}$. By
Theorem~\ref{expofthirdterm}, $v_i^{f_i(p^N)}$ will be trivial if
\[ N \geq a+\left\lfloor\frac{(k+1)-(s+1)}{p-1}\right\rfloor + \lfloor\log_p(s)\rfloor.\]
This must hold for $s=2,\ldots,k$. We set $s'=s-1$ and rewrite the above as:
\[ N \geq a + \left\lfloor\frac{(k-1)-s'}{p-1}\right\rfloor +
\lfloor\log_p(s'+1)\rfloor,\]
with $s'=1,\ldots,k-1$.
By Lemma~\ref{boundfors}, the largest value the right hand side takes is
$a + \left\lfloor\frac{k-1}{p-1}\right\rfloor$.
This proves that $[z^{p^N},y]=[z,y]^{p^N}$ if $N\geq a + \left\lfloor\frac{k-1}{p-1}\right\rfloor$.
The proof that $[z,y^{p^N}]=[z,y]^{p^N}$ is
essentially the same.
\end{proof}

\begin{remark} The results above also hold if the identity
  in~$(\mathrm{\ref{identityppower}})$ is replaced by an identity in
  which $p^i$ is placed on any of the three entries in $[z,y,z]$, and in
  any of the three entries in $[z,y,y]$; for instance,
\[ \relax[z^{p^i},y,z]=[z,y,y^{p^i}]=e,\quad\mbox{or}\quad
  [z,y^{p^i},z]=[z^{p^i},y,y]=e,\]
etc.
\end{remark}

\begin{lemma} Let $k\geq 1$ be a positive integer, $p$ a prime, and let
  $H$ be a nilpotent $p$-group of class~$k+1$. Suppose that
  $y_1,\ldots,y_r$ are elements of $H$ such that their images generate
  $H/Z(H)$; assume further that the order of $y_iZ(H)$
  in $H/Z(H)$ is~$p^{\alpha_i}$,
  with $1\leq \alpha_1\leq\cdots\leq \alpha_r$. Then $\alpha_r\leq
  \alpha_{r-1} + \left\lfloor\frac{k-1}{p-1}\right\rfloor$; that is,
\[ \smash{y^{p^{\alpha_{r-1}+\left\lfloor\frac{k-1}{p-1}\right\rfloor}}}\in Z(H).\]
\label{acentralelement}
\end{lemma}

\begin{proof}
Since $H$ is generated by $y_1,\ldots,y_r$ and central elements, it
is sufficient to prove that $[y_r^{p^{\alpha_{r-1}+\left\lfloor\frac{k-1}{p-1}\right\rfloor}},y_j]=e$ for
$j=1,\ldots,r-1$.  

Note that $y_j^{p^{\alpha_{r-1}}}$ is central for $j=1,\ldots,r-1$, so that
\[\forall i\geq \alpha_{r-1},\qquad
	  [y_r,y_j^{p^i},y_j]=[y_r,y_j^{p^i},y_r]=e.\]
Thus, we may apply Theorem~\ref{likenil2withpowers} to conclude that
\[ \relax
[y_r^{p^{\alpha_{r-1}+\left\lfloor\frac{k-1}{p-1}\right\rfloor}},y_j] =
[y_r,y_j^{p^{\alpha_{r-1}+\left\lfloor\frac{k-1}{p-1}\right\rfloor}}]=e,\]
proving our lemma.
\end{proof}

The necessary condition is now immediate:

\begin{theorem}[P.~Hall if $k<p$~\cite{hallpgroups}]
Let $G$ be a nilpotent $p$-group of class~$k$, with $k\geq 1$ and $p$
a prime. Furthermore, let $\{x_1,\ldots,x_r\}$ be a generating
set for~$G$ with $x_i$ of order $p^{\alpha_i}$, where
$\alpha_1\leq\alpha_2\leq\cdots\leq\alpha_r$. If $G$ is capable, then $r>1$ and
\[\alpha_{r}\leq\alpha_{r-1}+\left\lfloor\frac{k-1}{p-1}\right\rfloor.\]
\label{necessity}
\end{theorem}

\begin{proof}
Since a center-by-cyclic group is abelian, the necessity of $r>1$ is
clear. So assume that $G$ is capable, and $r>1$. Let $H$ be a
$p$-group of class $k+1$ such that $G\cong H/Z(H)$. Let
$y_1,\ldots,y_r$ be elements of $H$ that project onto $x_1,\ldots,x_r$,
respectively. 
Then Lemma~\ref{acentralelement} gives the condition on $\alpha_r$,
proving the theorem.
\end{proof}

\begin{remark} If $G$ is of small class then 
  $\lfloor\frac{k-1}{p-1}\rfloor=0$, so Theorem~\ref{necessity} says
  that $\alpha_r\leq \alpha_{r-1}$; therefore, when $k<p$ the
  necessary condition becomes ``$r>1$ and
  $\alpha_r=\alpha_{r-1}$.'' This is Hall's observation
  in~\cite{hallpgroups}, applied to a minimal generating set.
\end{remark}

The obvious question to ask is whether the inequality is tight. It is
clearly best possible if $k<p$. An easy case to consider for $k\geq p$
is $p=2$. In this case, the dihedral group of order $2^{k+1}$ is of
class~$k$, minimally generated by an element of order~$2$ and an
element of order $2^k$; and its central quotient is isomorphic to the
dihedral group of order $2^k$. Thus, the inequality is tight for all
$k$ when $p=2$. The case of $k\geq p>2$ is more difficult, but once
again there is a $2$-generator capable group where the equality holds,
with one generator of order~$p$. This is shown
in~\cite{inequality}. These considerations lead to the following
proposition.

\begin{prop}[\cite{inequality}] For every $k\geq 1$ and every prime~$p$ there exists
  a capable $p$-group of class~$k$, minimally generated by an element
  of order $p$ and an element of order
  $p^{1+\lfloor\frac{k-1}{p-1}\rfloor}$. In particular, the bound in
  Theorem~\ref{necessity} is best possible.
\end{prop}

\section{$k$-nilpotent products of cyclic $p$-groups,
  $k<p$\label{sec:centkprod}}

In this section we describe the center of a $k$-nilpotent product of
cyclic $p$-groups in the case $k\leq p$; as a result of this
description, we will prove the promised generalization of Baer's
Theorem for the small class case. The description is proven by
induction on~$k$, and we will need to do the case $k=2$ explicitly.

\begin{lemma} Let $p$ be a prime and $C_1,\ldots,C_r$ be cyclic $p$-groups
 generated by $x_1,\ldots,x_r$, respectively, with
  $p^{\alpha_i}$ the order of $x_i$ and $1\leq
  \alpha_1\leq\cdots\leq \alpha_r$. If $G=C_1\amalg^{{\germ N}_2}\cdots\amalg^{{\germ
  N}_2} C_r$, then \[Z(G) = \left\langle
  x_r^{p^{\alpha_{r-1}}}, G_2\right\rangle.\]
\label{centertwonilproduct}
\end{lemma}

\begin{proof} A theorem of T.~MacHenry~\cite{machenry} shows that the
  cartesian $[B,A]$ of $A\amalg^{{\germ N}_2} B$, with $A,B\in {\germ
  N}_2$, is isomorphic to $B^{\rm ab}\otimes A^{\rm ab}$, the
  isomorphism being the one 
  that sends $[b,a]$ to $\overline{b}\otimes\overline{a}$. Using the
  universal property of the coproduct, we map $G$ to
  $C_i\amalg^{{\germ N}_2} C_r$ to conclude that
  $C_i\cap
  Z(G) = \{e\}$ for $i=1,\ldots,r-1$. Looking at
  $C_{r-1}\amalg^{{\germ N}_2} C_r$ we also obtain that
  $C_r\cap Z(G) = \langle
  x_r^{p^{\alpha_{r-1}}}\rangle$, proving our claim.
\end{proof}

We will also need two observations on basic commutators:

\begin{lemma}
Let $F$ be the free group on $x_1,\ldots,x_r$. Let $[u,v]$ be a basic
commutator in $x_1,\ldots,x_r$, and assume that ${\rm wt}([u,v])=k\geq
2$.
\begin{itemize}
\item[(i)] If $v\leq x_r$, then $[u,v,x_r]$ is a basic commutator in
  $x_1,\ldots,x_r$.
\item[(ii)] If $v>x_r$, then
$[u,v,x_r] \equiv [v,x_r,u]^{-1}[u,x_r,v]  \pmod{F_{k+2}}$.
In addition, both $[v,x_r,u]$ and $[u,x_r,v]$ are basic
commutators in $x_1,\ldots,x_r$.
\end{itemize}
\label{jacobigeneral}
\end{lemma}

\begin{proof} Clause (i) follows from the definition of basic
  commutator, as does the claim in clause (ii) that $[v,x_r,u]$ and
  $[u,x_r,v]$ are basic commutators. The congruence in (ii) follows from
  Proposition~\ref{Widentities}(iv).
\end{proof}

\begin{lemma} Let $F$ be the free group on $x_1,\ldots,x_r$, and let
  $k\geq 2$. Let $c_1,\ldots,c_s$ be the basic commutators in
  $x_1,\ldots,x_r$ of weight exactly $k$ listed in ascending order,
  and write $c_i=[u_i,v_i]$.
  Let $\beta_1,\ldots,\beta_s$ be any integers, and let
  $g=c_1^{\beta_1}\cdots c_s^{\beta_s}$. For
  $i=1,\ldots,s$, let $d_i$ and $f_i$ be defined by:
\begin{eqnarray*}
d_i &=& \left\{\begin{array}{ll}
\relax[u_i,v_i,x_r]&\mbox{if $v_i\leq x_r$};\\
\relax[v_i,x_r,u_i]^{-1}&\mbox{if $v_i>x_r$.}
\end{array}\right.\\
f_i &=&\left\{\begin{array}{cl}
e&\mbox{if $v_i\leq x_r$}\\
\relax[u_i,x_r,v_i]&\mbox{if $v_i>x_r$.}
\end{array}\right.
\end{eqnarray*}
Then
\[
\relax [g,x_r]  \equiv \prod_{i=1}^s d_i^{\beta_i}f_i^{\beta_i}
\pmod{F_{k+2}}\]
and, except for removing the trivial terms with $f_i=e$,
this expression is in normal form for the abelian group
$F_{k+1}/F_{k+2}$.
\label{wecangodown}
\end{lemma}

\begin{proof} It is easy to verify that if $[u_i,v_i]\neq[u_j,v_j]$,
  then $d_i,f_i,d_j,f_j$ will be pairwise distinct, except perhaps in
  the case where $f_i=f_j=e$. Thus in the expression given for
  $[g,x_r]$, all terms are either powers of the identity or of pairwise
  distinct basic commutators. That the expression is indeed equal to
  $[g,x_r]$ follows from Proposition~\ref{Widentities}(ii).
\end{proof}

We are now ready to prove our result:

\begin{theorem} For a positive integer $k$ and a prime $p$ with $p\geq
  k$, let $C_1,\ldots,C_r$ be cyclic $p$-groups generated by
  $x_1,\ldots,x_r$ respectively, with $p^{\alpha_i}$ being the order
  of~$x_i$, and assume that $1\leq\alpha_1\leq \alpha_2\leq\cdots\leq
  \alpha_r$. If $G$ is the $k$-nilpotent product of the $C_i$, $G =
  C_1 \amalg^{\germ N_k}\cdots\amalg^{\germ N_k} C_r$, then $Z(G) =
  \langle x_r^{p^{\alpha_{r-1}}}, G_k\rangle$.
\label{centerknilprod}
\end{theorem}

\begin{proof} That the center contains the right hand side 
 follows from Lemma~\ref{acentralelement} and the properties of a
  $k$-nilpotent product. To prove the other inclusion, we proceed by
  induction on $k$. The result is trivially true for $k=1$, and
  Lemma~\ref{centertwonilproduct} gives the result for $k=2$.

Assume the result is true for the $(k-1)$-nilpotent product of the
$C_i$, $2<k\leq p$, and let $K=C_1 \amalg^{\germ N_{k-1}}\cdots
\amalg^{\germ N_{k-1}} C_r$; that is, $K=G/G_k$. 

By the induction hypothesis, we know that $Z(K)=\langle
x_r^{p^{\alpha_{r-1}}},K_{k-1}\rangle$, and the center of $G$
is contained in the pullback of this subgroup. So we have the
inclusions
$\langle x_r^{p^{\alpha_{r-1}}},G_{k}\rangle \subseteq Z(G)
\subseteq \langle x_r^{p^{\alpha_{r-1}}}, G_{k-1}\rangle$.
Let $g\in Z(G)$; multiplying by adequate elements of $G_k$ and an
adequate power of $x_r^{p^{\alpha_{r-1}}}$, we may assume that $g$ is
an element of $G_{k-1}$ which can be written in normal form as:
\[ g = \prod_{i=1}^{n} c_i^{\beta_i},\]
where $c_{1},\ldots,c_n$ are the basic commutators of weight
exactly $k-1$ in $x_1,\ldots,x_r$, and the $\beta_i$ are integers on the
relevant interval. If we prove that $g=e$, we will obtain our result.

Since $g\in Z(G)$, its commutator with $x_r$ is trivial. From
Proposition~\ref{Widentities}(ii) we have:
\[e  =  [g,x_r]
  =  \left[ \prod_{i=1}^n c_i^{\beta_i}, x_r\right]
 =  \prod_{i=1}^{n} [c_i,x_r]^{\beta_i}.\]

For each $i=1,\ldots,n$, write $c_i=[u_i,v_i]$, with $u_i,v_i$
basic commutators. Let $d_i$ and~$f_i$ be as in the statement of
Lemma~\ref{wecangodown}.  Then we have
\[e = [g,x_r] = \prod_{i=\ell}^n d_i^{\beta_i}f_i^{\beta_i}.\]
Except for some $f_i$ which are trivial, the precise ordering of the
remaining terms, and the exponents for the $d_i$ corresponding to
$v_i>x_r$, this is already in normal form. The ordering of the
nontrivial basic commutators is immaterial, since the $d_i$ and $f_j$
commute pairwise; and for those $d_i$ corresponding to $v_i>x_r$, we
simply add the corresponding power of $p$ to the exponent $-\beta_i$
to obtain an exponent in the correct range (see 
\cite{struikone}*{Theorem~3}).  The expression is then in normal form, so the only
way in which this product can be the trivial element is if
$\beta_i=0$ for $i=1,\ldots,n$, proving that $g=e$, and so the
theorem.
\end{proof}

We now have the desired generalization:

\begin{theorem}[Baer~\cite{baer} for $k=1$] 
For $k$ a positive integer and $p$ a prime with $p>k$, let
$C_1,\ldots,C_r$ be cyclic $p$-groups generated by
$x_1,\ldots,x_r$ respectively, where $p^{\alpha_i}$ is the order of
$x_i$ with $1\leq \alpha_1 \leq \alpha_2 \leq \cdots \leq
\alpha_r$. If $G$ is the $k$-nilpotent product of the $C_i$, $G = C_1
\amalg^{{\germ N}_k} \cdots \amalg^{{\germ N}_k} C_r$, then $G$ is
capable if and only if $r>1$ and $\alpha_{r-1}=\alpha_r$.
\label{capabilitynilkprod}  
\end{theorem}

\begin{proof} Necessity follows from Theorem~\ref{necessity}. For
  sufficieny, let 
\[ K = C_1\amalg^{{\germ N}_{k+1}}\cdots \amalg^{{\germ N}_{k+1}}
C_r.\]
Since $\alpha_{r-1}=\alpha_r$, we have from
Theorem~\ref{centerknilprod} that $Z(K)=K_{k+1}$, so
$K/Z(K)=K/K_{k+1}\cong G$, as desired.
\end{proof}

\begin{remark} The original statement of Baer's Theorem is not
  restricted to torsion groups, or even to finitely generated
  groups. The original result is:
\end{remark}

\noindent\textbf{Baer's Theorem} (Corollary to Existence Theorem
in~\cite{baer})\textbf{.}
\textit{
A direct sum $G$ of cyclic groups (written additively) is capable if
and only if it satisfies the following two conditions:
\begin{itemize}
\item[\textit{(i)}] If the rank of $G/G_{\rm tor}$ is 1, then the orders of the
  elements in $G_{\rm tor}$ are not bounded; and
\item[\textit{(ii)}] If $G=G_{\rm tor}$, and the rank of $(p^{i-1}G)_p/(p^iG)_p$
  is~$1$, then $G$ contains elements of order $p^{i+1}$, for all
  primes~$p$;
\end{itemize}
where $kG=\{kx\mid x\in G\}$, and $H_p = \{h \in H\mid ph=0\}$ for any
subgroup $H$ of~$G$.}

The reader may find it interesting to give a proof of Baer's Theorem
using our methods and the $2$-nilpotent product; a complete
description of the multiplication table for the $2$-nilpotent product
is given by Golovin in~\cite{metab}. The proof of Baer's Theorem is then
straightforward, and only complicated by the notation needed to
consider infinite direct sums.

\section{The case $k=p=2$\label{sec:kp2}}

In this section we consider the smallest case that is not covered by
our investigations so far: $k=p=2$. In this instance,
Theorem~\ref{necessity} gives the condition $\alpha_r\leq
\alpha_{r-1}+1$. We will prove that the condition is also sufficient
for the case of the $2$-nilpotent product of cyclic $2$-groups. 

As before, we start by examining the center of a $3$-nilpotent product
of cyclic $2$-groups. Such a product was considered in detail by
R.R.~Struik in~\cites{struikone,struiktwo}. To obtain uniqueness in
the normal form, we must replace the basic commutators $[z,y,z]$ and
$[z,y,y]$ with commutators $[z^2,y]$ and $[z,y^2]$, respectively, and
adjust the ranges of the exponents accordingly. The normal form result
we obtain with these changes is described in \cite{struikone}*{Theorem~4};
explicitly, it states that if $C_1,\ldots,C_r$ are cyclic groups
generated by $x_1,\ldots,x_r$ respectively, and if the order of $x_i$
is $2^{\alpha_i}$, $1\leq \alpha_1\leq \cdots \leq \alpha_r$, then
every element of the $3$-nilpotent product of the $C_r$ may be written
uniquely in normal form as:
\[ g  =   x_1^{\gamma_1} \cdots x_r^{\gamma_r} 
\!\!\!\!\prod_{1\leq i<j\leq r}
\!\!\!\![x_j,x_i]^{\gamma_{ji}} [x_j^2,x_i]^{\gamma_{jij}}[x_j,x_i^{2}]^{\gamma_{jii}}
\!\!\!\!\!\!\!\prod_{1\leq i<j<k\leq r}
\!\!\!\![x_j,x_i,x_k]^{\gamma_{jik}}[x_k,x_i,x_j]^{\gamma_{kij}},\]
%
where $\gamma_i, \gamma_{jik}$, and~$\gamma_{kij}$ are integers modulo
$2^{\alpha_i}$; $\gamma_{ji}$ is an integer modulo $2^{\alpha_i+1}$;
$\gamma_{jii}$ is an integer modulo $2^{\alpha_i - 1}$; and $\gamma_{jij}$ is an
integer modulo $2^{\alpha_i-1}$ if $\alpha_i=\alpha_j$, and modulo
$2^{\alpha_i}$ if $\alpha_i < \alpha_j$. 

Struik also provides multiplication formulas
in~\cite{struikone}*{pp.~453}, which may be verified by applying the
collection process; note however that our choice of basic commutators
differs slightly from hers, so the formulas in~\cite{struikone} do not
apply as is to our choice of normal forms. It is straightforward to do
the necessary conversions. In any case we will not need
the multiplication formulas here.

To simplify notation we set $\alpha=\alpha_{r-1}$. Let $G$ be as above;
we want to determine the center of~$G$. Clearly,
 $G_3\subset Z(G)$, and from Lemma~\ref{acentralelement} we know that
 $x_r^{2^{\alpha+1}}$ is central.
By
considering $G/G_3$, we obtain the $2$-nilpotent product of the $C_i$, and so we
know from Lemma~\ref{centertwonilproduct} that the center of $G$ satisfies:
\[ \left\langle x_r^{2^{\alpha+1}},G_3\right\rangle \subset
 Z(G) \subset \left\langle
 x_r^{2^{\alpha}},G_2\right\rangle.\]

We claim that in fact
\[ \left\langle x_r^{2^{\alpha+1}},G_3\right\rangle \subset
 Z(G) \subset \left\langle
 x_r^{2^{\alpha+1}},G_2\right\rangle.\]
This can be seen as follows: if $\alpha=\alpha_r$, then 
$x_r^{2^{\alpha+1}} = x_r^{2^{\alpha}} = e$,
and the claim is true. If, on the other hand, we have
$\alpha<\alpha_r$, then consider the commutator of $x_{r-1}$
with~$x_r^{2^{\alpha}}$. From~$(\ref{prodformthree})$
 we have:
\begin{eqnarray*}
\relax[x_r^{2^{\alpha}},x_{r-1}] &=&
      [x_r,x_{r-1}]^{2^{\alpha}}
      [x_r,x_{r-1},x_r]^{\binom{{2^{\alpha}}}{2}}\\
& = &
      [x_r,x_{r-1}]^{2^{\alpha}-2^{\alpha}}[x_r^2,x_{r-1}]^{2^{\alpha-1}}\\
& = & [x_r^2,x_{r-1}]^{2^{\alpha-1}};
\end{eqnarray*}
where we have used the identity
$[x_j,x_i,x_j]=[x_j,x_i]^{-2}[x_j^2,x_i]$, which may be obtained from
$(\ref{prodformthree})$ as well. Since $\alpha=\alpha_{r-1}<\alpha_r$, the
order of $[x_r^2,x_{r-1}]$ is $2^{\alpha}$, so this commutator
is not trivial. Therefore,
$x_r^{2^{\alpha_{r-1}}}$ is not central, but its square is, giving
once again the inclusions claimed.

Now let $g\in Z(G)$; multiplying by suitable powers of $x_r$ and of commutators of the form
$[x_j,x_i,x_k]$ and $[x_k,x_i,x_j]$ we may
assume that:
\[ g = \prod_{1\leq i<j\leq r}
   [x_j,x_i]^{\gamma_{ji}} [x_j^2,x_i]^{\gamma_{jij}} [x_j,x_i^2]^{\gamma_{jii}}.\]
We first take the commutator of $g$ with $x_r$; from
$(\ref{prodformone})$ and Proposition~\ref{Widentities}(ii) we have:
\begin{eqnarray*}
\relax [g,x_r] & = & \prod_{1\leq i<j\leq
  r}[x_j,x_i,x_r]^{\gamma_{ji}}[x_j^2,x_i,x_r]^{\gamma_{jij}}[x_j,x_i^2,x_r]^{\gamma_{jii}}\\
& = & \prod_{1\leq i< j\leq
  r}[x_j,x_i,x_r]^{\gamma_{ji}+2\gamma_{jij}+2\gamma_{jii}}.
\end{eqnarray*}
For $j<r$, we conclude that $\gamma_{ji}+2\gamma_{jij}+2\gamma_{jii}\equiv 0
\pmod{2^{\alpha_i}}$. If $j=r$ we again use the identity
$[x_r,x_i,x_r]=[x_r,x_i]^{-2}[x_r^2,x_i]$ and conclude that we must
have $-2(\gamma_{ri}+2\gamma_{rir}+2\gamma_{rii})\equiv 0 \pmod{2^{{\alpha_i}+1}}$ by
looking at the exponent of $[x_r,x_i]$ in the resulting expression. So
in any case we conclude that
\[  \gamma_{ji} + 2\gamma_{jij} + 2\gamma_{jii} \equiv 0 \pmod{2^{\alpha_i}};\qquad
1\leq i<j\leq r.\]
Conversely, it is easy to verify that if $g$ as above satisfies this
condition, then it will necessarily be
central in~$G$. We obtain:

\begin{theorem}
Let $C_1,\ldots,C_r$ be cyclic groups, generated by $x_1,\ldots,x_r$
respectively; assume that the order of $x_i$ is $2^{\alpha_i}$, and
that 
$1\leq \alpha_1\leq \alpha_2\leq \cdots \leq \alpha_r$. Let $G$ be the
3-nilpotent product of the $C_i$,
$G = C_1\amalg^{{\germ N}_3}\cdots \amalg^{{\germ N}_3} C_r$, and
write $\alpha=\alpha_{r-1}$.
Then $g\in Z(G)$ if and only if it can be written in normal form as:
\[ g = x_r^{\gamma_r} 
\!\!\!\!\!\prod_{1\leq i<j\leq r} 
\!\!\![x_j,x_i]^{\gamma_{ji}}[x_j,x_i^2]^{\gamma_{jii}}[x_j^2,x_i]^{\gamma_{jij}}
\!\!\!\!\!\prod_{1\leq i<j<k\leq r} 
\!\!\![x_j,x_i,x_k]^{\gamma_{jik}}[x_k,x_i,x_j]^{\gamma_{kij}},
\]
where $\gamma_r\equiv 0 \pmod{2^{\alpha+1}}$, and
$\rho_{ji}\equiv 0 \pmod{2^{\alpha_i}}$, where
$\rho_{ji}=\gamma_{ji} + 2\gamma_{jii} + 2\gamma_{jij}$.
That is
\[ Z(G) = \Biggl\langle \Bigl\{x_r^{2^{\alpha+1}},
   G_3\Bigr\}\;\bigcup\; \Bigl\{
   [x_j,x_i]^{2^{\alpha_i}}\,\Bigm|\, 1\leq i<j\leq r\Bigr\}\Biggr\rangle.\]
\label{centerofnil32groups}
\end{theorem}

With a description of the center, we can now easily derive the
characterization of the capable $2$-nilpotent products of cyclic
$2$-groups: 

\begin{theorem}
Let $C_1,\ldots,C_r$ be cyclic $2$-groups generated by
$x_1,\ldots,x_r$, respectively, where $2^{\alpha_i}$ is the order of
$x_i$, and assume that $1\leq \alpha_1\leq\cdots\leq \alpha_r$. If
$G$ is the $2$-nilpotent product of the $C_i$,
\[ G = C_1 \amalg^{{\germ N}_2} C_2 \amalg^{{\germ N}_2} \cdots
\amalg^{{\germ N}_2} C_r,\]
then $G$ is capable if and only if $r>1$ and $\alpha_r \leq \alpha_{r-1}+1$.
\label{pandkequal2}
\end{theorem}

\begin{proof} Necessity follows from Theorem~\ref{necessity}. For
  sufficiency, let $K$ be the $3$-nilpotent product of the $C_i$,
$K = C_1 \amalg^{{\germ N}_3} C_2 \amalg^{{\germ N}_3} \cdots
\amalg^{{\germ N}_3} C_r$.
Then the description of the center at the end of
Theorem~\ref{centerofnil32groups} makes it easy to verify that
$K/Z(K)\cong G$, so $G$ is capable.
\end{proof} 

The preceding theorem suggests the following question:

\begin{question} Is the $k$-nilpotent product of $r$ cyclic $p$-groups
  capable if and only if $r>1$ and 
\[\alpha_r\leq\alpha_{r-1}+\left\lfloor\frac{k-1}{p-1}\right\rfloor
?\]
\end{question}

As we have seen, the answer is yes when $k<p$, and when $k=p=2$.

\section{Some applications of our approach\label{sec:applic}}

As we noted in the introduction, one weakness of
Theorems~\ref{capabilitynilkprod} and~\ref{pandkequal2} is that
whereas the $1$-nilpotent product of cyclic groups covers all
finitely generated abelian groups, the case $k\geq2$ does not do the
same for the finitely generated nilpotent groups of
class~$k$. However, it is possible to use our results as a starting
point for discussing capability of other more general $p$-groups.

A recent result of Bacon and Kappe characterizes the capable
2-generated nilpotent $p$-groups of class two with $p$ an odd prime
using the nonabelian tensor square~\cite{baconkappe}. We can recover
their result using our techniques, and easily obtain a bit more.

In Theorem~2.4 of~\cite{baconkappenonab}, the authors present a
classification of the finite $2$-generator $p$-groups of class two,
$p$ an odd prime. With a view towards their calculations of the
nonabelian tensor square, the authors classify the groups into three
families. We will modify their classification and coalesce them into a
single presentation.

Let $G=\langle a,b\rangle$ be a finite nonabelian $2$-generator
$p$-group of class~$2$, $p$ an odd prime. Then $G$ is isomorphic to
the group presented by:
\begin{equation}
\left\langle a,b\,\Biggm|\,
\begin{array}{rcl}
a^{p^{\alpha}}=b^{p^{\beta}}=[b,a]^{p^{\gamma}} &=&e,\\
\relax[a,b,a]=[a,b,b] &=& e,\\
a^{p^{\alpha+\sigma-\gamma}} [b,a]^{p^{\sigma}}& = & e.
\end{array} \right\rangle
\label{eq:presentation}
\end{equation}
where $\alpha+\sigma\geq 2\gamma$, $\beta\geq\gamma\geq 1$,
$\gamma\geq \sigma\geq 0$,
and if $\sigma=\gamma$,
then $\alpha\geq\beta$.  Under these restrictions, the choice is
uniquely determined. 

If $\sigma=\gamma$, we obtain the groups in Bacon and Kappe's first
family, which one might call the ``coproduct type'' groups (they are
obtained from the nilpotent product $\langle a\rangle\amalg^{{\germ
N}_2}\langle b\rangle$ by taking the quotient modulo a power of $[a,b]$,
and are
the coproduct in a suitably chosen subvariety of~${\germ N}_2$). If
$\sigma=0$, we obtain the split meta-cyclic groups, which are the
second family in~\cite{baconkappenonab}. The cases $0<\sigma<\gamma$
correspond to their third family.

The result which appears in~\cite{baconkappe} is that a
$2$-generated group with presentation as in~(\ref{eq:presentation})
with $\sigma=0$ or $\sigma=\gamma$ is capable if and only if
$\alpha=\beta$. The condition is also both necessary and sufficient
for the remaining case with $0<\sigma<\gamma$ (in this
case,~\cite{baconkappe} contains an error which the
authors are in the process of correcting~\cite{kappepers}).

Thus, in the case of $2$-generated $p$-groups of class~two, $p$ an odd
prime, Baer's condition is both necessary and sufficient, just as for
finite abelian groups.

Although we could prove the result for all three families at once, we
will divide the argument in two, to prove slightly more for
the case where $\sigma=\gamma$.

Let $p$ be an odd prime and $C_1,\ldots,C_r$ be cyclic $p$-groups
generated by $x_1,\ldots,x_r$, respectively, where $p^{\alpha_i}$ is
the order of $x_1$ with $1\leq \alpha_1\leq\alpha_2\leq\cdots\leq
\alpha_r$. Let $K$ be the $2$-nilpotent product of the $C_i$,
\[ K = C_1 \amalg^{{\germ N}_2}\cdots \amalg^{{\germ N}_2} C_r.\]
For each pair $j,i$, let $\beta_{ji}$ be a positive integer less than
or equal to~$\alpha_{i}$, and set $N=\langle
[x_j,x_i]^{p^{\beta_{ji}}}\rangle$. Since $N$ is central, we have
$N\triangleleft K$. Let $G=K/N$.

\begin{theorem} The group $G$, as defined in the previous paragraph,
is capable if and only if $r>1$ and $\alpha_{r-1}=\alpha_r$.
\end{theorem}

\begin{proof} Necessity follows from Theorem~\ref{necessity}. For
 sufficiency, let $H$ be the $3$-nilpotent product of the $C_i$,
$H = C_1\amalg^{{\germ N}_3}\cdots \amalg^{{\germ N}_3} C_r$.
Clearly $H$ will not do, so we need to take a quotient of $H$ so that,
in the resulting group, $[x_j,x_i]^{\beta_{ji}}$ is central. To that
end, we let $M$ be the subgroup of $H$ generated by all elements of the form
$[x_j,x_i,x_k]^{\beta_{ji}}$ with $1\leq i<j\leq r$, and $k$ arbitrary.
In terms of the basic commutators, this means the elements
$[x_j,x_i,x_k]^{\beta_{ji}}$ for $k\geq i$, and the elements
\[[x_i,x_k,x_j]^{-\beta_{ji}}[x_j,x_k,x_i]^{\beta_{ji}}\]
for $k<i$.  It is easy to verify that all elements of~$H/M$ have a
normal form as in Theorem~3 of~\cite{struikone}, and that the
multiplication of these elements uses the same formulas as those in
$H$, except that the exponents of the basic commutators of weight~$3$
are now taken modulo the adequate $\beta_{ji}$ instead of the old moduli.

Proceeding now as in the proof of Theorem~\ref{centerknilprod}, one proves
that the center of $H/M$ is generated by the third term
of the lower central series, the image of $x_r^{p^{\alpha_{r-1}}}$,
and those of the elements
$[x_j,x_i]^{p^{\beta_{ji}}}$; therefore by taking the central quotient, 
we obtain the group~$G$, as desired.
\end{proof}

The first part of \cite{baconkappe}*{Corollary~4.3} corresponds to
the case $r=2$, $x_1=a$, $x_2=b$, and $\beta_{21}=\gamma$.
Now we consider the case with $\sigma<\gamma$.

\begin{theorem}[cf.\ \cite{baconkappe}*{Corollary~4.3}]
 Let $p$ be an odd prime, and let $G$ be a group presented by
 $(\ref{eq:presentation})$ with $0\leq\sigma<\gamma$. Then $G$ is capable
 if and only if $\alpha=\beta$.
\end{theorem}

\begin{proof} Necessity once again follows from Theorem~\ref{necessity}; so we
  only need to prove sufficiency. The idea is to construct the
``obvious witness'' to the capability, by starting with the
$3$-nilpotent product of two cyclic groups of order $p^{\alpha}$,
generated by~$x$ and~$y$; then for every relation $r$ in the
presentation of~$G$, we will take the quotient modulo the normal
closure of $\langle [r,x],[r,y]\rangle$, thus making $r$ central in
the quotient. Then we just need to make sure that in the resulting
group, the map $x\mapsto a$, $y\mapsto b$ will yield the desired
isomorphism between the central quotient and~$G$.  In essence, what we
are doing is constructing a ``generalized extension of~$G$'' which can
be used to determine the capability
of~$G$. See~\cite{beyltappe}*{Theorem~III.3.9}.

So let $H = \langle x\rangle \amalg^{{\germ N}_3} \langle y\rangle$,
where $x$ and $y$ are both of order $p^{\alpha}$. 

First, we want to make sure that $[y,x]^{p^{\gamma}}$ is central.
To that end, we let $N=\langle [y,x,x]^{p^{\gamma}},[y,x,y]^{p^{\gamma}}\rangle$,
and consider $K=H/N$.

The next step is to ensure that
$x^{p^{\alpha+\sigma-\gamma}}[y,x]^{p^{\sigma}}$ is central. So first
we consider
\begin{eqnarray*}
[x^{p^{\alpha+\sigma-\gamma}}[y,x]^{p^{\sigma}},x] & = &
[x^{p^{\alpha+\sigma-\gamma}},x] [x^{p^{\alpha+\sigma-\gamma}},x,[y,x]^{p^{\sigma}}]
[[y,x]^{p^{\sigma}},x]\\
& = & [y,x,x]^{p^{\sigma}}
\end{eqnarray*}
(using (\ref{prodformone}) and Proposition~\ref{Widentities}). So we
let $L=K/\langle [y,x,x]^{p^{\sigma}}\rangle$. 

So far, the only difference in the normal form and
multiplication tables for $L$ and for $H$ is the order of
$[y,x,x]$ and $[y,x,y]$.

The final quotient we need to take is to ensure that
$x^{p^{\alpha+\sigma-\gamma}}[y,x]^{p^{\sigma}}$ also commutes with
$y$. To that end, we consider:
\begin{eqnarray*}
[x^{p^{\alpha+\sigma-\gamma}}[y,x]^{p^{\sigma}},y] & = &
 [x^{p^{\alpha+\sigma-\gamma}},y] [x^{p^{\alpha+\sigma-\gamma}},y,[y,x]^{p^{\sigma}}]
[[y,x]^{p^{\sigma}},y]\\
& = & [x^{p^{\alpha+\sigma-\gamma}},y] [y,x,y]^{p^{\sigma}}\\
& = & [y,x]^{-p^{\alpha+\sigma-\gamma}}
 [y,x,x]^{-\binom{p^{\alpha+\sigma-\gamma}}{2}} [y,x,y]^{p^{\sigma}}.
\end{eqnarray*}
The last equality uses $(\ref{prodformthree})$.
Note, however, that
the conditions on $\alpha$, $\sigma$, and~$\gamma$ imply that
$\alpha+\sigma-\gamma>\sigma$, so we can
simplify the expression above to:
\[ [x^{p^{\alpha+\sigma-\gamma}}[y,x]^{p^{\sigma}},y] 
= [y,x]^{-p^{\alpha+\sigma-\gamma}}[y,x,y]^{p^{\sigma}}.\]
Let $N=\langle
[y,x]^{p^{\alpha+\sigma-\gamma}}[y,x,y]^{-p^{\sigma}}\rangle$. 
It is easy to verify that
for all $g\in L$, both $[g,x]$ and $[g,y]$ lie in $N$ if and only if 
\[  g \in \left\langle x^{p^{\alpha+\sigma-\gamma}}[y,x]^{p^{\sigma}},
    [y,x]^{p^{\gamma}}, [y,x,x], [y,x,y]\right\rangle,\] (for example,
using the explicit multiplication formulas in~\cite{struikone},
suitably modified to fit our choice of basic commutators) from which
we deduce that if $M=L/N$, then $M/Z(M)\cong G$, as desired.
\end{proof}


Let $p$ be an odd prime, and let $G$ be a noncyclic $2$-generator finite
$p$-group of class at most two. If $G=\langle a,b\rangle$, then at
least one of $a$ and $b$ must have the same exponent as the exponent
of the element of highest order of the group, since $G$ is regular;
say $a$, with order $p^{\mu}$. Among all elements $b'\in G$ such that
$G=\langle a,b'\rangle$, pick the one of smallest order and replace
$b$ with $b'$; say the order is $p^{\nu}$. Then $\mu$ and $\nu$ are
invariants of~$G$; in fact, they are the two largest type invariants
of the group (see~\cite{hallcontrib}*{\S4}). Arguing as in
\cite{baconkappenonab}*{Lemma 2.3}, one shows that $\langle
a\rangle\cap\langle b'\rangle = \{e\}$, and then following the
argument in~\cite{baconkappenonab}*{Th.~2.4}, it is easy to
verify that $\mu$ and $\nu$ are equal to $\alpha$ and $\beta$ (in some
order) in the presentation given in (\ref{eq:presentation}). So we obtain:

\begin{corollary}
Let $p$ be an odd prime, and let $G$ be a noncyclic $2$-generator
finite $p$-group of class at most two, and let $\alpha$ and $\beta$ be
the two largest type invariants of~$G$. Then $G$ is capable if and
only if $\alpha=\beta$.
\label{cor:refsugg}
\end{corollary}

At this point, an obvious question to ask is whether the necessary
condition of Theorem~\ref{necessity} will also prove to be necessary
for the case of nilpotent groups of class~two, as it was for 
abelian groups (at least, for odd prime and a suitable choice of minimal generating
set). Unfortunately, the answer to that question is no.

Recall that a $p$-group $G$ is \textit{extra-special} if and only if
$G'=Z(G)$, $|G'|=p$, and $G^{\rm ab}$ is of exponent~$p$.  The theorem
of Beyl, Felgner, and Schmid in~\cite{beyl} mentioned in the
introduction states that an extra-special $p$-group is capable if and
only if it is dihedral of order order~$8$, or of order $p^3$ and
exponent $p$, with $p>2$. So the extra-special $p$-group $G$ of order
$p^5$ generated by $x_1, x_2, x_3, x_4$,  and satisfying%
\ $[x_3,x_1]=[x_3,x_2]=[x_4,x_1]$, $[x_4,x_2]=[x_4,x_3]=[x_2,x_1]=e$,
$x_1^p=x_2^p=x_3^p=x_4^p=e$, and $G_2\subset Z(G)$, 
is not capable, and is minimally generated by four elements of exponent
$p$. Thus, the necessary condition is not sufficient in general for
groups in~${\germ N}_2$.
An independent proof that $G$ is not capable together with some other applications
using our methods for groups of exponent~$p$ and class~$2$ 
appears in~\cite{capablep}.

For now, we note that the nilpotent product can be used to produce a
natural ``candidate for witness'' to the capability of a given
finite nilpotent $p$-group $G$. Let $G$ be a finite nilpotent group of
class~$k>0$, minimally generated by elements $x_1,\ldots,x_r$. Let
$w_1(x_1,\ldots,x_r),\ldots,w_n(x_1,\ldots,x_r)$
be words in
$x_1,\ldots,x_r$ that give a presentation for $G$, that is:
\[ G = \Bigl\langle x_1,\ldots,x_r \,\Bigm|\,
w_1(x_1,\ldots,x_r),\ldots,w_n(x_1,\ldots,x_r)\Bigr\rangle.\]
Let $\langle y_1\rangle,\ldots,\langle y_r\rangle$ be infinite cyclic
groups, let 
$K=\langle y_1\rangle\amalg^{{\germ N}_{k+1}}\cdots \amalg^{{\germ
    N}_{k+1}} \langle y_r\rangle$,
and let $N = \Bigl\langle
w_1(y_1,\ldots,y_r),\ldots,w_n(y_1,\ldots,y_r)\Bigr\rangle^K$,
that is, the normal closure of the subgroup generated by the words
evaluated at $y_1,\ldots,y_r$. Finally, let $M=[N,K]$.

\begin{theorem}
The group $G$, as defined in the previous paragraph, is capable if and
only if
\[ G \cong (K/M)\bigm/ Z(K/M),\]
where $K$ and $M$ are also as defined in the previous paragraph.
\label{th:specialwitness}
\end{theorem}

\begin{proof} We only need to prove the ``only if'' part. Note that the
  map that sends $y_1,\ldots,y_r$ to $x_1,\ldots,x_r$ respectively
  gives a well-defined map from $K$ to $G$, and that $K/\langle
  N,K_{k+1}\rangle\cong G$; since $M\subseteq N$, 
\[(K/M)/\langle NM,K_{k+1}M\rangle\cong G,\]
and $\langle NM, K_{k+1}M\rangle$ is central in $K/M$. Therefore,
  $(K/M)/Z(K/M)$ is a quotient of $G$. We want to show that it is
  actually isomorphic to~$G$.

Since $G$ is capable, there is a group $H$ such that $H/Z(H)\cong
G$. Let $h_1,\ldots,h_r$ be elements of $H$ that map to
$x_1,\ldots,x_r$, respectively. Replacing $H$ by $\langle
h_1,\ldots,h_r\rangle$ if necessary, we may assume that
$h_1,\ldots,h_r$ generate $H$. Since $G$ is of class~$k$, $H$ is of
class~$k+1$, and therefore there exists a unique surjective map
from~$K$ (the relatively free ${\germ N}_{k+1}$ group of rank $r$) to
$H$, mapping $y_1,\ldots,y_r$ to $h_1,\ldots,h_r$, respectively. We
must have that $w_i(h_1,\ldots,h_r)\in Z(H)$ for $i=1,\ldots,n$, so
the map from $K$ to $H$ factors through $K/M$. Since
$Z(K/M)$ maps into $Z(H)$, the induced mapping $K/M\to H\to G$
factors through $(K/M)/Z(K/M)$. Therefore, we have that $(K/M)/Z(K/M)$
has $G$ as a quotient. Since it is also isomorphic to a quotient of
the finite group~$G$, we must have $(K/M)/Z(K/M)\cong
G$, as claimed. \end{proof}

\section*{Acknowledgements}
I became interested in the question of capability after attending 
talks by L.C. Kappe and M.~Bacon at the AMS Meeting in Baltimore. They
were kind enough to give me a copy of their preprint
of~\cite{baconkappe}, and L.C.~Kappe very patiently answered my
questions over e-mail. I thank them both very much for their time and
help.

In addition, several people were kind enough to help me during the
preparation of this paper. Bill Dubuque, Ed Hook, Robert Israel,
and~Joseph Silverman helped with some of the calculations relating to
binomial coefficients. George Bergman offered helpful comments.
Avinoam Mann brought several known results to my attention.  Derek
Holt, Mike Newman, and Vasily Bludov from the Group Pub Forum verified
some of my hand-calculations using MAGMA. Finally, Robert Morse used
GAP to verify some of my constructions, answered many of my questions,
and also provided some counterexamples to guesses of mine. I thank all
of them very much for their time and their help.

Finally, I am very grateful to the referee, who offered many
suggestions that greatly improved readability, as well as pointing out a gap in
the original proof of Lemma~\ref{boundfors} and suggesting
the statement of Corollary~\ref{cor:refsugg}.

\section*{Errata. Added 4/10/2006. Updated 5/15/2006.}

The last clause of Lemma~4.2(ii) as stated above is unfortunately false.
For example, if the weight of $[v,x_r]$ is smaller than the
weight of $u$, then $[v,x_r,u]$ cannot be a basic commutator; in
addition, if $u$ is of the form $[c_1,c_2]$ with $c_i$ a basic
commutator and the weight of $c_2$ is greater than $1$, then
$[u,x_r]$ is not a basic commutators, hence neither
is $[u,x_r,v]$. Analogous problems arise if $v$ is of this form.

The problem can be resolved by continuing the process, using work of
Ward and of T.C.~Hurley; a proof of Theorem~4.4 follows from those
arguments. I am preparing a second part which will contain the
correction of the problem, as well as extending the results to $k=p$
for arbitrary prime~$p$.

\end{document}